\documentclass{svproc}
\usepackage{url}

\usepackage{scalerel}
\usepackage{hyperref}
\usepackage{type1cm}        
\usepackage{makeidx}         
\usepackage{graphicx}        

\usepackage{multicol}        
\usepackage[bottom]{footmisc}

\usepackage{newtxtext}       %
\usepackage[varvw]{newtxmath}       

\usepackage{wrapfig}

\usepackage{bm}
\usepackage{siunitx}

\usepackage{graphicx, subcaption} 
\usepackage{amsmath, hyperref, mathtools, color, caption}
\usepackage{array}
\usepackage[capposition=bottom]{floatrow}
\usepackage{algpseudocode}
\usepackage{ulem}
\usepackage{multirow, makecell, booktabs, verbatim}

\usepackage{algorithm}

\DeclareSymbolFont{bbold}{U}{bbold}{m}{n}
\DeclareSymbolFontAlphabet{\mathbbold}{bbold}

\begin{document}

\mainmatter              
\title{WoS-NN: an Effective Stochastic Solver for Elliptic PDEs with Machine Learning}
\titlerunning{WoSNN: Stochastic Solver for PDEs with Machine Learning}  
%
\author{Silei Song \and Arash Fahim \and Michael Mascagni}
\authorrunning{S. Song, A. Fahim, M. Mascagni} 
%
\tocauthor{Silei Song, Arash Fahim, and Michael Mascagni}
\institute{Florida State University, Tallahassee, FL 32306, USA, 
\email{ss19cu@fsu.edu}}
\maketitle              
\begin{abstract}
Solving elliptic partial differential equations (PDEs) is a fundamental step in various scientific and engineering studies. As a classic stochastic solver, the Walk-on-Spheres (WoS) method is a well-established and efficient algorithm that provides accurate local estimates for PDEs. In this paper, by integrating machine learning techniques with WoS and space discretization approaches, we develop a novel stochastic solver, WoS-NN. This new method solves elliptic problems with Dirichlet boundary conditions, facilitating precise and rapid global solutions and gradient approximations. The method inherits excellent characteristics from the original WoS method, such as being meshless and robust to irregular regions. By integrating neural networks, WoS-NN also gives instant local predictions after training without re-sampling, which is especially suitable for intense requests on a static region. A typical experimental result demonstrates that the proposed WoS-NN method provides accurate field estimations, reducing errors by around $75\%$ while using only $8\%$ of path samples compared to the conventional WoS method, which saves abundant computational time and resource consumption.

\end{abstract}

\section{Introduction}
In the modern development of science and engineering, partial differential equations (PDEs) are the basis of various natural phenomena and industrial applications. Solving elliptic PDEs such as Laplace and Poisson equations enhances our understanding of natural processes and drives advancements across multiple industries \cite{SalsaPDE}. Classical deterministic numerical solvers for PDEs include the finite element method, finite difference method, domain decomposition method, Galerkin method, etc. \cite{numeric}. These techniques are based on the differential or integral transformations of the original PDEs and typically require meshing and calculations across the entire domain, which is time-consuming. On the other hand, non-deterministic methods such as stochastic solvers simulate true solutions via sampling in the domain \cite{Sabelfeld1991}. They provide sufficiently accurate local estimations without requiring meshing, which makes them advantageous in practical scenarios with lousy boundary conditions and discontinuity situations.

As a stochastic solver for elliptic PDEs, our WoS-NN method combines the classic WoS method with modern machine-learning techniques. In our study, PDEs are reformulated through their stochastic representations,  where WoS is employed as a random path generator and a spatial discretization tool. The WoS method offers fast, discretized random walks to boundaries and avoids the time factor in the approximation process. On the other hand, neural networks (NNs) are used to approximate the PDE solutions and gradients along the WoS paths. After training with the WoS sample paths, the networks can estimate solutions and gradients of the target PDE for any new local points in $O(1)$ time.

In section 2, a background on elliptic PDE solvers is provided, including the original WoS method and recent machine-learning approaches in PDE solving. In Section 3, we present the primary processes of the WoS-NN method, including the mathematical background and our neural network design. Section 4 compares WoS-NN with other stochastic solvers on solving different PDEs under diverse conditions. Experiments demonstrate significant advantages of WoS-NN over the original WoS method. Finally, Section 5 summarizes our work and outlines some potential future improvements.

\subsection{Contribution}

In this paper, we propose a fast stochastic solver, WoS-NN, of elliptic PDEs to give field estimations for both solution and gradient. Our method adopts the WoS algorithm to discretize stochastic processes and constructs a novel recurrent network structure for the approximation. As a Monte Carlo method, WoS-NN inherits WoS's advantages: no meshing, fast convergence, flexibility, and capability to handle intricate regions. Moreover, it offers accurate global approximations for both solutions and gradients with limited random walks, giving immediate estimations for new local solutions. The Recurrent-Neural-Network-like (RNN-like) model is compatible with varying-length sampling paths, providing smoothed field estimations and additional gradient approximations for noise reduction. As a result, the WoS-NN solver achieves accurate instant estimations for PDE solutions and can be widely applied in areas like geometry processing, physics simulation, electrical engineering, biomolecular modeling, etc. \cite{Rohan2,Capacitance1991,Simonov2004,Bossy2010}


\section{Background}

 Our paper focuses on second-order elliptic PDEs with Dirichlet boundary conditions. In this paper, we specifically focus on Laplace equations and Poisson equations:
    \begin{align}
    \text{Laplace:}\quad\Delta u = 0 \quad \text{on} \quad \Omega\text{, and }\quad
    u = g \quad\text{on} \quad \partial \Omega\label{songeqn:Laplace}\\
    \text{Poisson:}\quad\Delta u = f \quad \text{on} \quad \Omega \text{, and }\quad
    u = g \quad \text{on} \quad \partial \Omega.\label{songeqn:Poisson}
    \end{align}

Solving these two equations with our WoS-NN method is crucial, as they are the foundation of solving other elliptic PDEs. Extensive experiments and detailed analyses of these two equations have demonstrated the significant advantages of our method over the original Walk on Spheres method and other stochastic approaches.

\subsection{Walk-on-Spheres Method}
The original Walk-on-Spheres method (WoS) was first proposed in 1956 \cite{Muller1956}. This method is grounded in the Mean Value Property and Kakutani's principle \cite{Kakutani1944} for the Laplace equation \eqref{songeqn:Laplace}. Kakutani's principle provides an expected representation of the original PDE as 
    $u(x_0) = \mathbf{E}[g(W_\tau)]$
, where $W_\tau \in \partial\Omega$ is the first exit point on the region boundary, reached by a random walk (Brownian motion) $W_t$ starting from $x_0$, with hitting time at $\tau$. 

\begin{figure}[H]
    \begin{subfigure}{0.45\textwidth}
    \centering
        \includegraphics[width=0.85\textwidth]{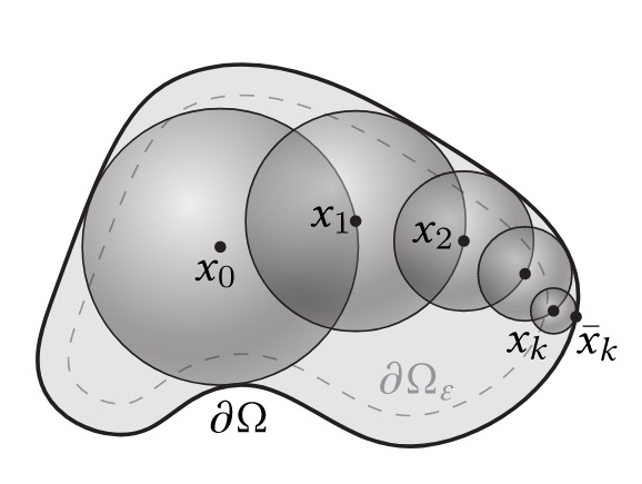}
        \caption{Walk on Spheres repeatedly moves to a random point on the largest sphere centered at the current point $x_k$, until it gets within an $\epsilon$ distance to the boundary. }
    \end{subfigure}
    \hfill
    \begin{subfigure}{0.45\textwidth}
    \centering
        \includegraphics[width=0.89\textwidth]{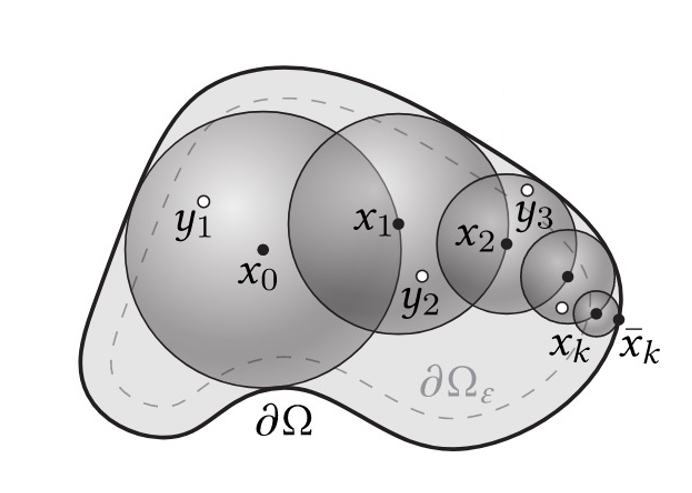}
        \caption{An additional random sample $y_{k + 1}$ is required inside each ball $B(x_k)$ to query the source term $f$ while solving Poisson problems.}
    \end{subfigure}
    \caption{The conventional WoS and samplings inside balls, from \cite{Rohan1}}\label{songfig:WoS}
\end{figure}

The WoS method leverages the Kakutani principle on a sphere by recursively sampling the exit point of  Brownian motion from internal spheres. The hitting time of Brownian motion starting at the center of a sphere is uniformly distributed on the surface of the sphere. More precisely, $x_i$, the exit point of Brownian motion from the $i$-th sphere, is the center of the $(i+1)$-th sphere, as shown in the left panel in Fig.\ref{songfig:WoS}. Thus, the Brownian motion $W_t$ is cut by spheres along the path. To approach the boundary $\partial \Omega$ as fast as possible, new centers $x_i$ are sampled from the largest possible spheres uniformly as the next WoS steps until the boundary is met. The uniform sample on the final ball is an evaluation of the Kakutani principle. An $\epsilon$-neighborhood inside the boundary has been adopted as a threshold determining if the new center reaches the boundary or not. If the new sample ($x_k$ in Fig.\ref{songfig:WoS} Left) is $\epsilon$ distance to the boundary, the boundary value of the closest point on the boundary ($\bar{x_k}$ in Fig.\ref{songfig:WoS} Left) will be adopted as a new estimation of the local solution at $x_0$. As a result, the estimation given by WoS can be represented by $\hat{u}(x_0) = \frac{1}{N} \sum_{i=0}^N g(\bar{x}^i_{k_i})$, where $N$ is the total amount of WoS sample paths, and $\bar{x}^i_{k_i}$ is the reached boundary point in $k_i$ step for the $i$-th path.
\vspace{-7mm}
The time complexity of WoS sampling each path is $O(\log(\frac{1}{\epsilon}))$, which is the average number of steps in a WoS trajectory with the $\epsilon$ capture region, proven in Theorem 8.1 \cite{Sabelfeld2013}. WoS paths converge rapidly to the boundary, regardless of the region's size and shape. The approximation accuracy of WoS depends on both the $\epsilon$-shell and the total sample path amount $N$, which can be adjusted according to demand. A smaller $\epsilon$ provides more precise boundary values to the Kakutani expectation, and a larger sample size is a common way to reduce the approximation variance for Monte Carlo methods.

One crucial improvement on the classical WoS method is to solve Poisson \eqref{songeqn:Poisson} and screened Poisson equations \cite{Rohan1}. Screened Poisson equations are also referred to as linearized Poisson-Boltzmann equations, which have extra source and absorption terms compared to the Laplace problems. The stochastic representations for screened Poisson equations are easily obtainable \cite{dynkin1965}. To estimate the source and absorption changes inside each WoS sphere, an additional in-ball sampling for each WoS step is needed \cite{Elepov1969}, as shown in the right panel in Fig. \ref{songfig:WoS}. Methods in \cite{Bossy2010,Simonov2007} offer the WoS extension on solving screened Poisson equations via implementing a killing term based on the in-ball sampling technique.

\subsection{Machine Learning in PDE Solving}

The rapid development of machine learning techniques has brought a new perspective to PDE solving, especially for high-dimensional, nonlinear problems. Current machine learning applications in PDE solving include combining machine learning techniques with other existing PDE solvers, such as the Deep Galerkin Method \cite{Sirignano2017} and the deep BSDE method \cite{Weinan2017,Beck2019}. The deep BSDE model solves high-dimensional non-linear parabolic PDEs by constructing deep neural networks on the time discretizations of Backward Stochastic Differential Equations (BSDEs). Another popular machine learning PDE solver is the Physics Informed Neural Networks \cite{RAISSI2019}, which incorporates physical laws directly into the training loss and collaborates with the operator learning. Compared to PINN, WoS-NN achieves efficient and stable convergence on elliptic PDE solving by incorporating the physical laws with the network structures, giving fast and intensive responses to new requests and new environments.

Our WoS-NN method applies machine learning approximation to the spatial discretization of elliptic PDEs. Moreover, rather than meshing the whole region, we employ WoS to discretize the random paths, which avoids temporal factors and allows flexible approximations of solutions and gradients. Other WoS and neural network combinations include using neural networks to cache the Monte Carlo estimation results and smoothing the approximation \cite{Siggraph1,Siggraph2}. In this method, neural networks are utilized as a variance reduction technique and can be widely applied to general Monte Carlo methods. Another incorporation of neural networks and WoS is NWoS \cite{Nam2024}, which is the improved realization of training neural networks with WoS's results on high-dimensional Poisson problems. Compared to these methods, WoS-NN innovatively adopts a recurrent network model and realizes synchronous gradient approximations with Ito's Lemma and second-order Taylor extension. The integrated training of both solutions and solution gradients enables the networks to realize more precise global optimization. 


\section{Our Method}
A typical stochastic solver consists of two steps: transforming the PDE into an SDE and then approximating the SDE's solution via stochastic process samplings. In WoS-NN, the first step is accomplished by Ito's Lemma, while the second step contains WoS spatially discretizing the region and a neural network to approximate solutions and gradients. In this section, we will introduce the stochastic transformations for different PDEs with the WoS method as the spatial discretization tool and our network structuring for each WoS step's approximation.

\subsection{Stochastic Representation of PDEs and WoS Discretization}

One theory adopted to transfer PDEs to stochastic differential equations is Ito's Lemma \cite{Ito}. It states how an Ito process performs on the 2nd-order Taylor expansion for any twice differentiable continuous function $u(x, t)$, where the differential of $u$ is composed of a time differential $dt$ and a differential on Brownian motion $dW_t$. When applying Ito's Lemma, a stochastic representation for a continuous and twice differentiable function $u$ is available as: 
\begin{equation}
    du(W_t) = \frac{1}{2}\Delta u(W_t)dt + \nabla u(W_t)dW_t \quad on \quad \Omega, \label{songeqn:ItoLemma}
\end{equation}
where $W_t$ is the simplest Ito process as a Brownian motion in the region: $W_t \in \Omega, 0 \leq t \leq \tau$. 

For the solution $u$ of Laplace equation \eqref{songeqn:Laplace}, we have $\Delta u(W_t) = 0, \forall W_t \in \Omega$. Thus, the stochastic representation of the Laplace equation's solution $u$ is as follows:
\begin{equation}
    du(W_t) = \nabla u(W_t)dW_t, \; \; \forall t \in [0, \tau], \; W_t \in \Omega.
\end{equation}

Same as in the Kakutani principle, we generate Brownian motions from internal local points $W_0 \in \Omega$ to the boundary as $W_\tau \in \partial\Omega$. Thus, for any fine enough time discretization $t_0 = 0 < t_1 < ... < t_N = \tau$ with $\tau$ as the stopping time of $W_t$ (the moment $W_t$ first hit $\partial \Omega$ at $t =t_N \; s.t. \; t_N = \tau$), the following approximation can be achieved:
\begin{equation}
    u(W_{t_{i+1}}) \approx u(W_{t_i}) + \nabla u(W_{t_i})\Delta W_{t_i}, \;\; \Delta W_{t_i} = W_{t_{i+1}} - W_{t_i}.
\end{equation}
And the reached boundary condition $u(W_\tau) = u(W_{t_n})$ can be estimated backwardly as:
\begin{equation}
    u(W_\tau) \approx u(W_0) + \scaleobj{.8}{\sum_{i = 0}^{n-1}} \nabla u(W_{t_i})\Delta W_{t_i}, \;\; \Delta W_{t_i} = W_{t_{i+1}} - W_{t_i} \label{songeqn:Lap-est}
\end{equation}

Since the Laplacian operator zeroes out the temporal factor brought by the Ito process, only the discretized Brownian motion $W_{\tau_i}$ is left to be approximated. As a result, a spatial discretization on $W_t$ is achieved on this region as equation \eqref{songeqn:Lap-est} instead of the time discretization. We construct neural networks to approximate local solutions and local gradients along the random process as $\hat{u} = YNN(\cdot, \theta)$, $\nabla \hat{u} = ZNN(\cdot, \theta)$ to obtain the following estimation:
\begin{equation}\label{songeqn:Lap-NN}
    \hat{u}(W_\tau) = \hat{u}(W_{t_N})= YNN(W_0, \theta) + \scaleobj{.8}{\sum_{i = 0}^{N-1}} ZNN(W_{t_i}, \theta)\Delta W_{t_i},
\end{equation}
where $\hat{u}$ is the estimation for the true solution $u$ based on this discretization and neural network approximations.

The loss function to train networks $YNN$ and $ZNN$ with sampled random paths $W_t \in \Omega$ is established on the differences between the estimated boundary condition $\hat{u}(W_\tau)$ and the true boundary condition as:
\begin{equation}\label{songeqn:Lap-loss}
    loss = \mathbf{E}[|\hat{u}(W_\tau) - g(W_\tau)|^2],
\end{equation}
which trains the random variables $\theta$ in equation \ref{songeqn:Lap-NN}, and optimizes network $YNN$ and $ZNN$ to fit the PDE's solution geometry and solution's gradients separately.

A significant advantage of applying this network design to elliptic PDEs' spatial discretization is that only two global approximations $YNN, ZNN$ are needed for the static system's simulation without the temporal factor, rather than having sub-networks to capture the changing situations for each time step \cite{Weinan2017}. Moreover, while the network structuring is no longer restricted by the time steps, our training data can be regarded as sampled paths with dynamic lengths under different discretization principles. Thus, the sampling process is much more flexible. 

\subsubsection{WoS for Spatial Discretization}

Instead of sampling the whole path before the discretization or meshing the region, the WoS-NN method adopts WoS as the spatial discretization tool to cut random paths. Using spheres to locate the sample steps along path $W_t$ is cheap and effective, and the gaps $\Delta W_{\tau_i}$ are always flexible and controllable by adjusting the largest spheres allowed. Within this discretization, each cut of the random path will be a WoS step $W_{t_i}=x_i$, and the gap between steps is the sphere radius $\Delta W_{t_i} = R_i=x_{i+1} - x_i$. The approximation equation \eqref{songeqn:Lap-NN} now is:

\begin{equation}\label{songeqn:YNN-ZNN-Laplace}
    \hat{u}(W_\tau) = YNN(x_0, \theta) + \sum_{i = 0}^{n-1} ZNN(x_i, \theta)R_i.
\end{equation}

\begin{figure}
\vspace{-4mm}
    \centering
        \includegraphics[height=43mm]{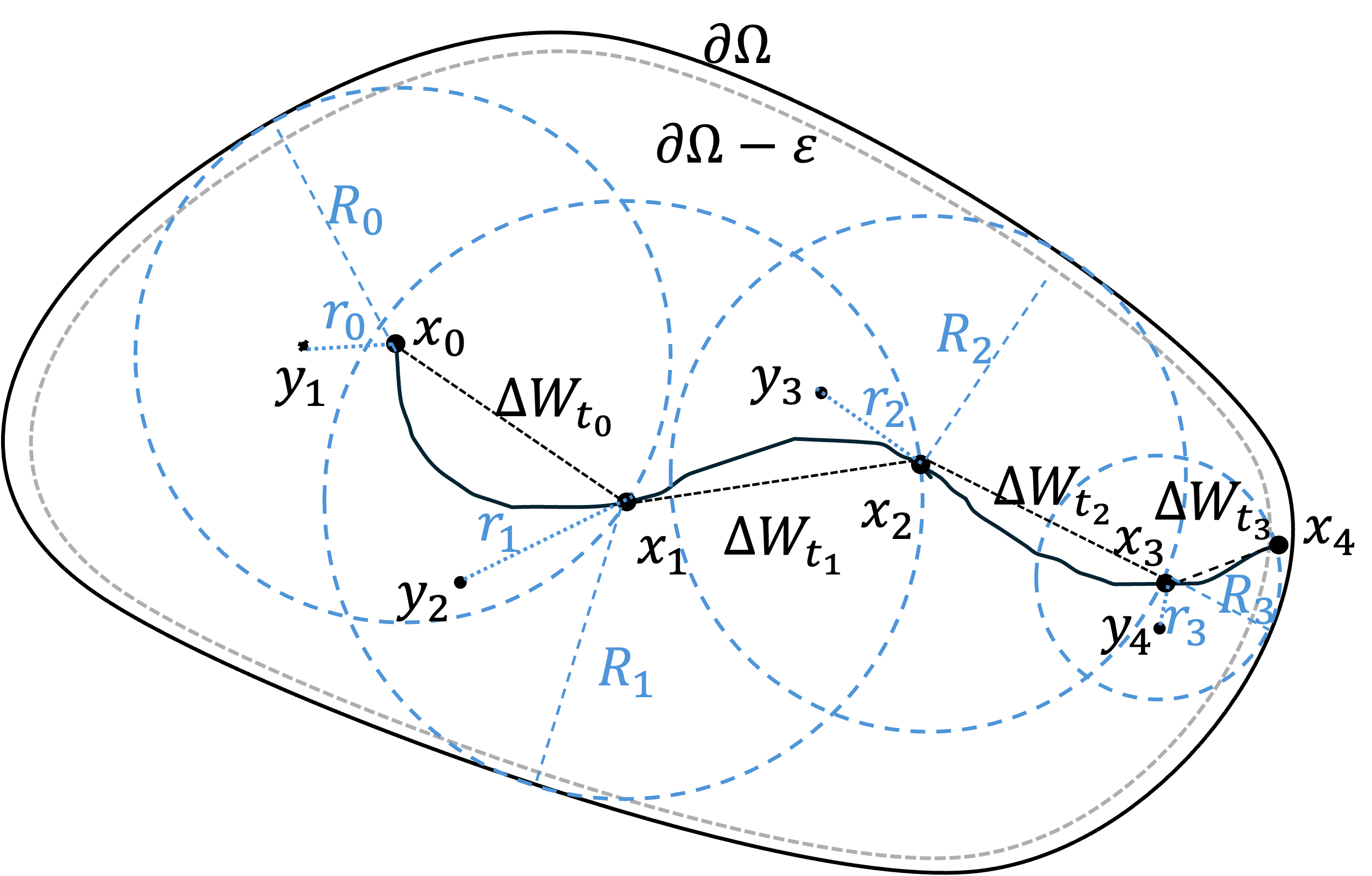}
    \caption{WoS for spatial discretization with each step cutting the random path at $W_{t_i} = x_i$.}
    \label{songfig:WoS_Poisson}
\end{figure}

As in Fig. \ref{songfig:WoS_Poisson}, the Brownian motion is split by spheres of a WoS process from the starting place $x_0$ to the reached boundary point $x_4$. The gradient changes along the path are accumulated inside each sphere as $\nabla u(W_{t_i})\Delta W_{t_i} \approx ZNN(x_i, \theta) R_i$. Thus, WoS-NN skips the meshing and significantly accelerates the approximation by utilizing WoS to take samples on spheres. 

Errors of our approximation can be introduced from the $\epsilon$-shell, the sample path amount, the maximum step amount allowed in a path, and the discretization step size. A proper maximum step limitation is necessary for different regions to balance between the errors brought by large walking steps and the accumulated gradient errors from long and winding paths.

\subsubsection{In-ball Samplings for Poisson Equation}
A stochastic representation of the Poisson equation \eqref{songeqn:Poisson} is derived from Ito's Lemma similarly as:
\begin{equation}
    du(W_t) = \frac{1}{2}\Delta u(W_t)dt + \nabla u(W_t)dW_t = \frac{1}{2}f(W_t)dt + \nabla u(W_t)dW_t.
\end{equation}
In contrast to Laplace equations, the time factor remains in Poisson's stochastic form, accumulating source terms along the Brownian motion. The integral representation for this stochastic equation is:
\begin{equation}\label{songeqn:Ito-Poisson}
    u(W_{t_{i+1}}) = \frac{1}{2}\scaleobj{.8}{\int_{t_i}^{t_{i+1}}}f(W_t)dt + \scaleobj{.8}{\int_{t_i}^{t_{i+1}}}\nabla u(W_t)dW_t.
\end{equation}
Thus, a joint sampling should be considered for both the time and Brownian motion in order to discretize the above equation \eqref{songeqn:Ito-Poisson} \cite{Deaconu2013}, to keep both $\Delta t_i = t_{i+1} - t_i$ and $\Delta W_{t_i} = W_{t_{i+1}} - W_{t_i}$ minor at the same time. 

To avoid the expensive joined sampling, paper \cite{Elepov1969} gave a way to approximate the source-term contribution $\int_{t_i}^{t_{i+1}}f(W_t)dt$ in any time interval $\Delta t_i$, with an extra sample $y_{i+1}$. They showed that:
\begin{equation}\label{songeqn:in-ball}
    \mathbf{E}\Big[\scaleobj{.8}{\int_{t_i}^{t_{i+1}}} f(W_t)dt\Big] = \mathbf{E}[f(y_{i+1})], 
\end{equation}
where $y_{i+1}$ is distributed proportional to the Green's function of the region $|y_{i+1}-x_i|\le R_i$. While being applied in WoS, the interval from $t_i$ to $t_{i+1}$ denotes the duration of Brownian motion $W_t$ within the i-th sphere, where $x_i = W_{t_i}$ is the sphere center. $F_{R_i}$ is a coefficient dependent solely on the radius of the i-th sphere, given as $R_i = \Delta W_{t_i}$. The in-ball sample $y_{i+1} \in B(x_i)$ is sampled with a probability $\mathbf{P}_{G(x_i)}$ according to Green’s function's density within the ball. 
Thus, with the source term integral estimated by in-ball samplings as in \eqref{songeqn:in-ball}, the WoS-NN can be extended to solve Poisson equations with the following estimation:
    \begin{equation}
        u(W_{t_{i+1}}) \approx u(W_{t_i}) + \frac{1}{2}f(y_{i+1})F_i + \nabla u(W_{t_i})R_i,
    \end{equation}
where $y_{i+1}$ is the extra in-ball sample for each step. Still, neural networks are used to approximate local solutions and local gradients as $u = YNN(\cdot, \theta)$ and $\nabla u = ZNN(\cdot, \theta)$ to have the following approximation: 
    \begin{equation}\label{songeqn:YNN-ZNN-Poisson}
        \hat{u}(W_\tau) = YNN(W_0, \theta) + \sum_{i=0}^{n-1}ZNN(W_i, \theta)R_i + \frac{1}{2}\sum_{i=0}^{n-1}f(y_i)F_i.
    \end{equation}
The loss function here is the same as the loss equation \eqref{songeqn:Lap-loss}, as for the moment Brownian motion hits the boundary. 

\subsection{Network Design}
\begin{figure}[h]

    \begin{subfigure}{0.45\textwidth}
    \centering
        \includegraphics[width=0.85\textwidth]{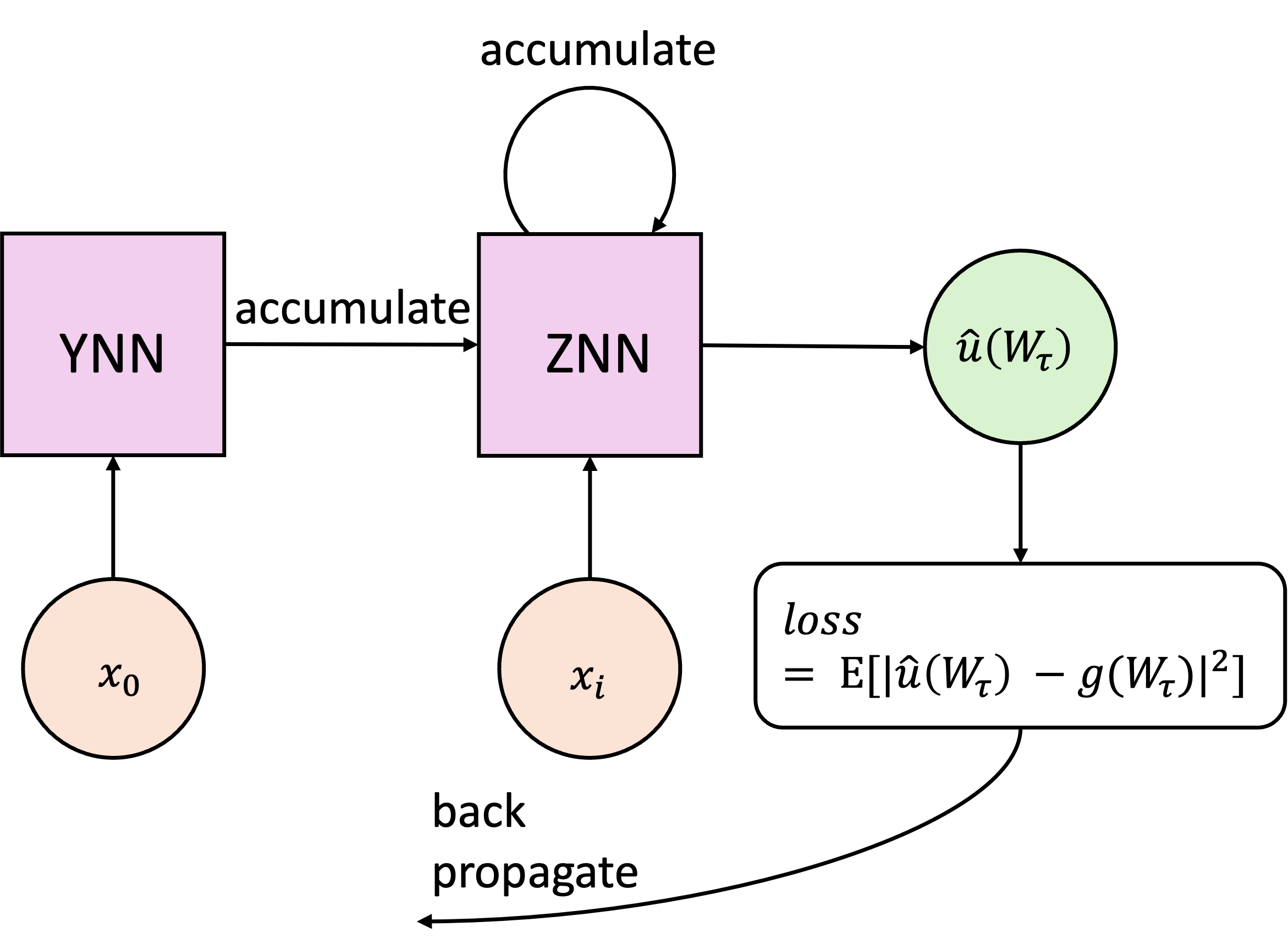}
        \caption{The initial NN design with two independent sub-networks. YNN estimates the solution at the starting place. ZNN estimates the gradient along paths.}
    \end{subfigure}
    \hfill
    \begin{subfigure}{0.45\textwidth}
    \centering
        \includegraphics[width=0.85\textwidth]{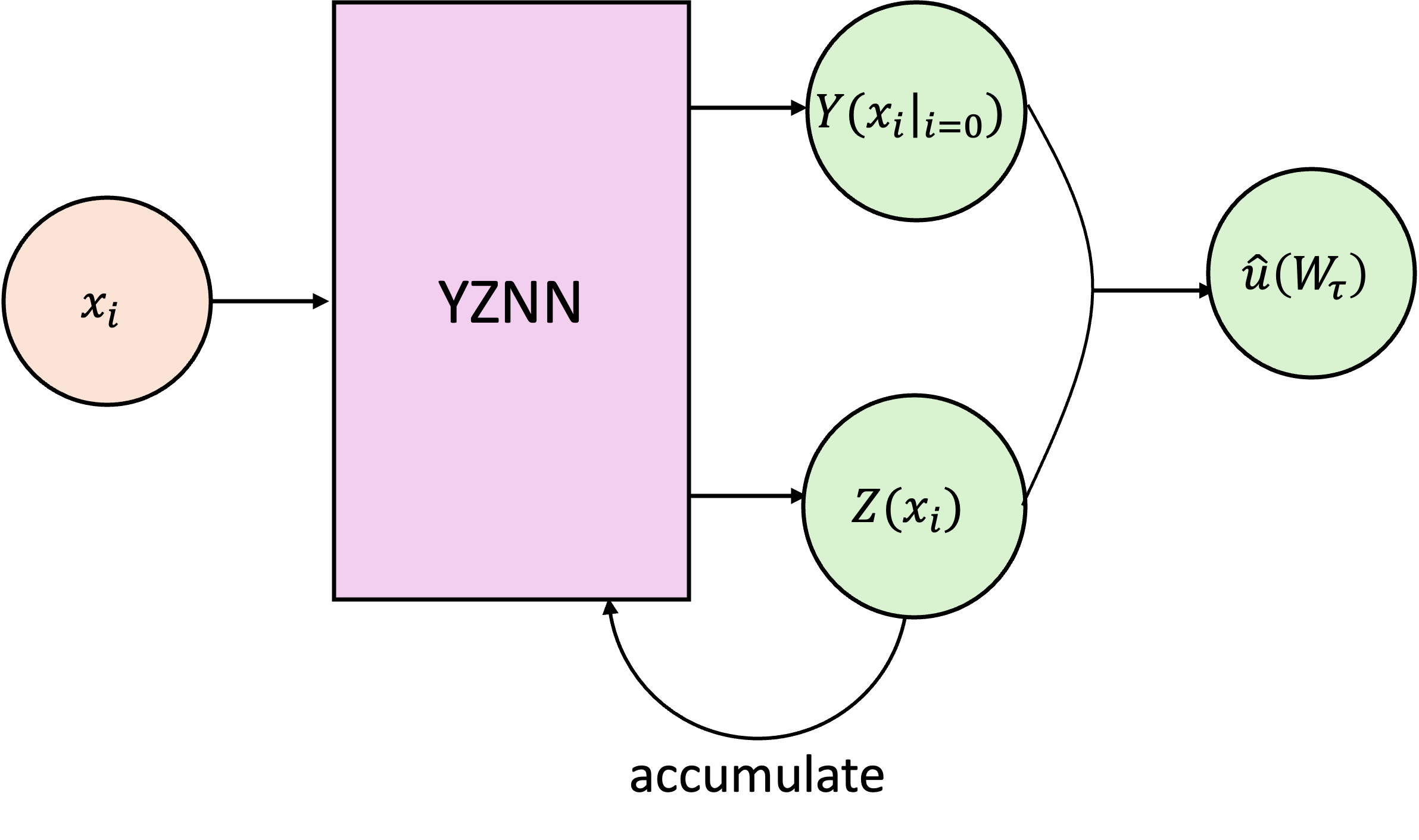}
        \caption{The improved NN design with one integrated network. Network YZNN simultaneously gives local approximations for both the solution and the gradient.}
    \end{subfigure}
    \caption{Two different WoS-NN neural network structures.}\label{songfig:subNN}
\end{figure}

As in equations \eqref{songeqn:YNN-ZNN-Laplace} and \eqref{songeqn:YNN-ZNN-Poisson}, we initially attempted to use two independent sub-networks to approximate local solutions and gradients separately. Our initial network structure is as in the left panel of Fig. \ref{songfig:subNN}. However, in this design, each random path trains network $YNN(\cdot, \theta)$ only once at the beginning but train $ZNN(\cdot, \theta)$ multiple times during every WoS step. This unbalanced training will force the optimization to focus on one of the two sub-networks, especially for winding paths with dynamic lengths. To fix the problem, we integrate two sub-networks into one, ensuring that the training of one approximation would affect the other. The improved neural network design is as in Fig. \ref{songfig:subNN} right panel. The new network $YZNN(\cdot, \theta)$ gives $n+1$-dimensional predictions for each local point, including a $1$-dimensional solution approximation and an $n$-dimensional gradient approximation ($\Omega \subseteq \mathbf{R}^n$), simultaneously.  

\begin{figure}[ht!]
    \centering
        \includegraphics[height = 63mm]{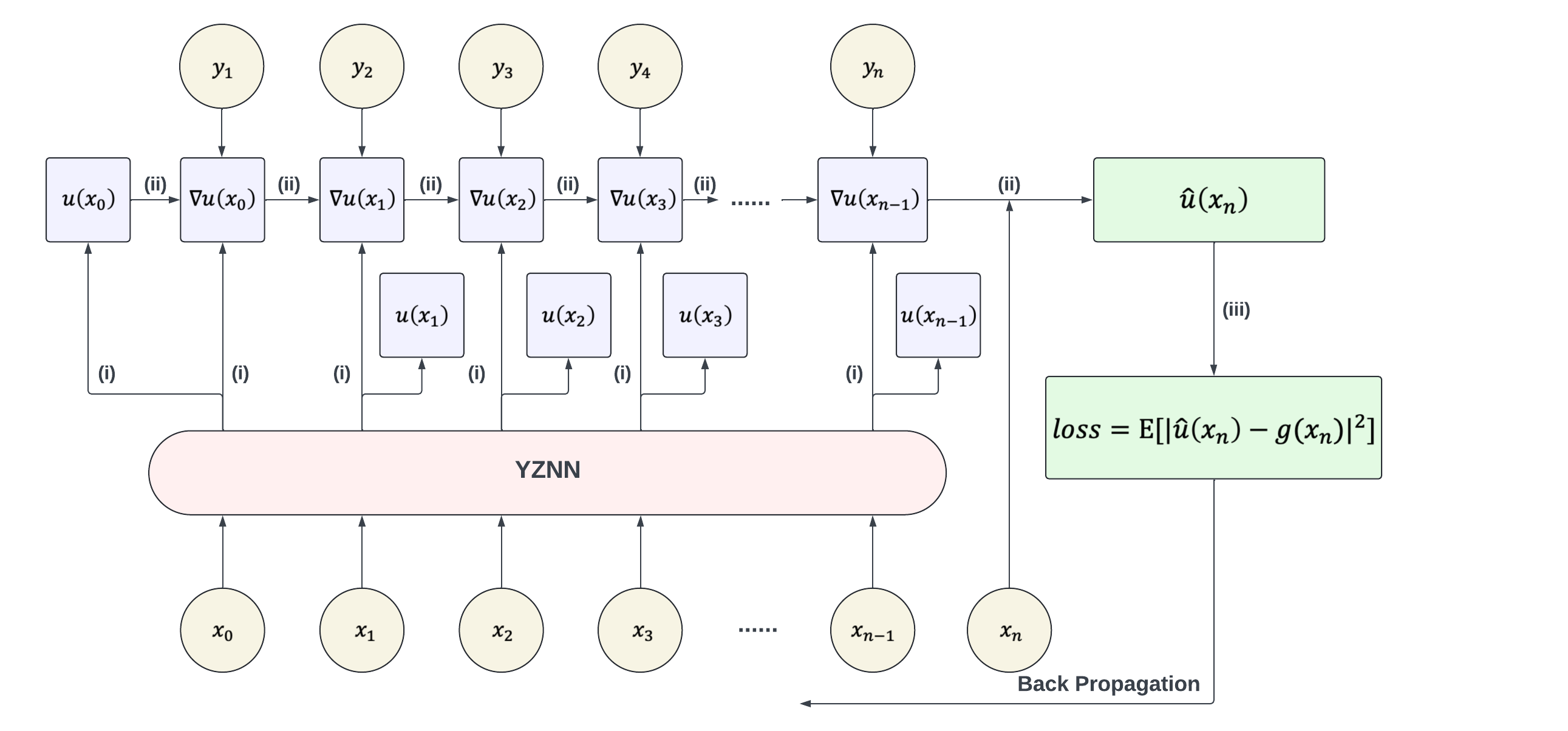}
    \caption{Structure of the loss function.}
    \label{songfig:loss}
\end{figure}

Fig.\ref{songfig:loss} is a comprehensive structure of our final network design, which gives an intuitive impression of the whole structure. The input training data is the sampled WoS discretized paths $x_i \in \Omega, i = 0, ..., n$. In Fig.\ref{songfig:loss}, steps (i) are the fitting steps that output solution and gradient approximations for any local $x_i$. Steps (ii) are simple linear calculations, imitating the gradient and source term changes along the path. When the path ends, a boundary estimation is provided as the output of the whole network, which is evaluated in the loss function (iii) for training. The $YZNN$ network acts as a recurrent network module whose feed-forward delivery can be repeated arbitrarily multiple times depending on the step amount. After training with enough samples, the network $YZNN$ can give precise solution approximations and gradient approximations at any position $x \in \Omega$ in $O(1)$ time. A general workflow of WoS-NN in Pseudocode form is as Alg. \ref{songalg:WoSNN}.

\begin{algorithm}[h]
\caption{WoS-NN for 2-dimensional Poisson equations}\label{songalg:WoSNN}
\begin{algorithmic}
\State For region $\Omega$, given Dirichlet boundary condition $g(x)$ on $\partial\Omega$ and source term $f(x)$
\State Neural network $YZNN_\theta$ with initial parameter $\theta$, batch size $m$ training for $T$ epochs
\State 1. Generate starting places $Pts$ and sample WoS paths
\For{$x_0 \in Pts$}
    \State Start point $X \gets x_0$, path $p$, source $s \gets 0$
    \While{$X \notin \partial \Omega_\epsilon$}
        \State Append $X$ to current path $p$
        \State Draw largest neighborhood $B(X, r)$, $r = distance (X, \partial \Omega)$
        \State Sample $Y \sim \mathbf{P}_{G(X)}(B(X, r))$, $s \gets  s-|B(X, r)|f(Y)G(X, Y)$
        \State Sample $X' \sim \mathcal{U}(\partial B(X, r))$  
        \State $X \gets X'$ 
    \EndWhile
    \State Find $\bar{X} \in \partial \Omega_\epsilon$ closest to $X$ 
    \State Record local approximation $\hat{u}^N(x_0) = g(\bar{X}) - s$ to $boundaries$ and path $p$ to $paths$
\EndFor
\State 2. Vectorize $paths$ by extending shorter paths to a uniform length $maxLen$. Shorter paths wait at the boundary.
\State 3. Training $YZNN$ with $paths$ and $boundaries$
\For{Every mini-batch in every epoch, path matrix $M_{m * maxLen * 2} \in paths$, corresponding boundary $B_b \in boundaries$}
    \State Local solution estimation for starting places as $loss \leftarrow YZNN(M[:, 0, :])[:, 0]$
    \State $loss \leftarrow YZNN(M[:, 0, :])[:, 0] + \sum_{i = 1}^{maxLen}YZNN(M[:, i, :])[:, 1:2] \cdot (M[:, i+1, :] - M[:, i, :])$
    \State $loss \leftarrow |loss - B_b|^2$
    \State $\theta \leftarrow SGD(loss)$
\EndFor
\end{algorithmic}
\end{algorithm}

In general, our WoS-NN method follows three main steps: the stochastic representation of PDEs, discretization of stochastic processes, and neural network approximations. The novel contributions of our method are as follows:
\begin{itemize}
        \item  We use a single neural network $YZNN$ for the static global estimation instead of having multiple sub-networks for each time step, since the Laplace operator eliminated the time factor from the stochastic system. As a RNN-like network, $YZNN$ can process paths with varying lengths without the limitation of time steps.
        \item For discretization, the WoS is adopted as a sampling tool to avoid meshing. Also, by using in-ball samplings in the WoS method, we avoid introducing the time factor back into the Poisson system.
        \item The neural network training spreads WoS local solutions to the neighborhood and offers a global approximation. As a field estimator, the solution and gradient approximations can be obtained anywhere from the trained network in constant time. 
        \item The gradient approximation performs as an additional noise reduction along each WoS path compared to the traditional WoS method, which smooths the solution estimation globally. Thus, WoS-NN is expected to give better results when the same sample paths are used.
        \item The integrated $YZNN$ network reveals the inner correlation between PDE solutions and their gradients. Instead of letting sub-networks capture this relationship automatically, combining the two sub-networks forces the network to pay attention to this inherent relationship among samples.
    \end{itemize}


\section{Experiments}

In this section, WoS-NN is tested with different parameters and factors under various PDE conditions. The experiments were executed in 2-dimensional and 3-dimensional regions for Laplace and Poisson equations. The Finite Difference method (FDM) and the original WoS method were tested for comparison on the same problems. We also trained a simple feed-forward neural network using WoS results directly as a reference. The experiment results indicate that WoS-NN performs better in various tests and examples with fewer sample paths. All experiments were done using Python 3.9, with a MacBook Pro machine holding an Apple M3 Pro microchip (11-core CPU, 14-core GPU, 18 GB RAM).

Our experiments were set up on 2-dimensional region $\Omega = [-1, 1]^2$ and 3-dimensional region $\Omega = [-1, 1]^3$. Models will be constructed and tested on an evenly distributed grid on the region as $Pts = \{(0.02 * i, 0.02 * j), i,j = \{-49, -48, ... 48, 49\}\} \in [-1, 1]^2$ and $Pts = \{(0.02 * i, 0.02 * j, 0.02 * k), i, j, k = \{-49, -48, ... 48, 49\}\}$, as of density $0.02 * 0.02$ or $0.02 * 0.02 * 0.02$. 

For 2-dimensional WoS-NN experiments, WoS paths were sampled with a maximum step of 20 and $\epsilon-$shell of 0.001. Each WoS-NN starting point only generated one random path, which was then filtered by the maximum step limitation. $YZNN$ is a fully connected feed-forward neural network with three hidden layers $(32, 64, 128)$. The network is trained with a $3*10^{-4}$ learning rate, ReLU activation function, and the Adam optimizer. Each training took 50 epochs to converge. For comparison, the original WoS method was run directly on the test set $Pts$ with a maximum step of 20, $\epsilon-$shell 0.001. Also, a WoS-driven neural network was trained directly with the WoS outputs on $Pts$. The WoS-driven network has the same structure as $YZNN$ and was trained for 50 epochs, with a $10^{-4}$ learning rate, batch size 256, ReLU activation function, and Adam optimizer, and was tested on the same target set $Pts$. Errors for different models were calculated as the average Euler distances between the estimated geometry and the actual solution on testing points $Pts$.

\subsection{2-dimensional Laplace Equations}
The 2-dimensional Laplace equation we've tested is: 
\begin{equation*}
        \Delta u = 0 \quad \text{on} \quad \Omega=(-1, 1)^2\text{, and  }\quad
        u = xy \quad \text{on} \quad \partial\Omega.
\end{equation*}
Table \ref{songtab:Laplace} compares representative experiment results for different models, where for each model, we take the fewest sample paths while the model converges.
\begin{table}[ht!]
    \centering
    \begin{tabular}{ |>{\centering\arraybackslash}m{2em}|>{\centering\arraybackslash}m{9em}|>{\centering\arraybackslash}m{6em}|>{\centering\arraybackslash}m{5em}|>{\centering\arraybackslash}m{5em}|>{\centering\arraybackslash}m{5em}| } 
 \hline
  &Method & Total Valid Path Amount & Execution Time \footnote{Execution time for NN-related models is denoted as sampling, training, and testing time} & Training Loss & Average Error on $Pts$ 
   \\ 
 \hline
 1 & Original WoS & 477459 & 69.3253 &  & 0.0377\\
 \hline
 2 & WoS-driven NN & $477459$  \footnote{We derived 9801 pairs of training data from 477459 WoS sample paths. As a result, the training path amount and the sampling time of experiment 2 are the same as in experiment 1.}& $69.3253$ $1.5884$ $2.5866$ & 0.0025 & 0.009475\\
 \hline
 3 & WoS-NN with uniform starting places & 38400 & $8.4687$ $48.8989$ 
$2.399$ & 0.02358 & \makecell{$0.008926$ \footnote{Solution approximation errors for WoS-NN.}\\[-0.8ex] \cmidrule(l{-4pt}r{-4pt}){1-1} $0.01497$ \footnote{Gradient approximation errors for WoS-NN.}}\\
 \hline
    \end{tabular}
    \caption{Experiment data for different models solving the Laplace equation.}
    \label{songtab:Laplace}
    
\end{table}

In Experiment 3, the WoS-NN was tested with WoS paths sampled from 40000 random starting places. The $YZNN$ network was then trained by these paths with a batch size of 2048 for 50 epochs. For comparison, experiment 1 ran the original WoS method on $Pts$, with at most 50 walks sampled per target local. The WoS-driven network was trained with the WoS output for $Pts$ in experiment 2, with batch size 256, for 50 epochs. 

Table \ref{songtab:Laplace} shows that our WoS-NN method provides competitive approximations with only $8\%$ sampling walks compared to both the original WoS method and the WoS-driven network. The training time for WoS-NN was longer due to its recurrent network structure. The WoS-NN network is trained on all points of each sample path. Different sample paths may vary in the number of steps. However, the training of WoS-NN is a one-time process, while the original WoS method resamples paths for every new local. Once the training is finished, WoS-NN can estimate both the solution and the gradient for any local place within $O(1)$ time. In the WoS-driven NN experiment, there are only 9801 pairs of $(x, \hat{u}_{WoS}(x))$ training data to train a simple feed-forward neural network on the starting points of each path, which leads to a short training time. However, since the training data of WoS-driven NN came directly from the classical WoS result, it requires the same sampling time as in the WoS experiment, which dominates the execution time of experiment 2.

\begin{figure}[ht!]
    \begin{subfigure}{0.44\textwidth}
    \centering
        \includegraphics[width=1\textwidth]{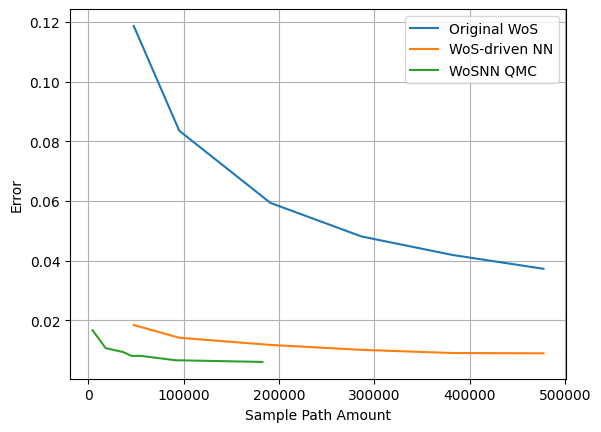}
    \caption{Mean error versus sample path amount}
    \end{subfigure}
    \quad
    \begin{subfigure}{0.42\textwidth}
    \centering
        \includegraphics[width=1\textwidth]{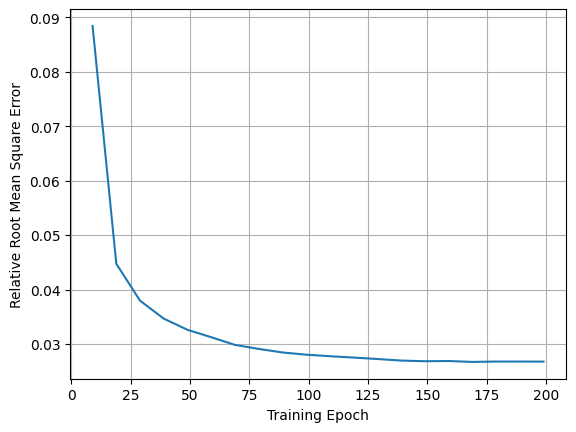}
    \caption{RRMSE versus the training epochs}
    \end{subfigure}
    \caption{Analytical experiment results for the Laplace experiments.}\label{songfig:line-chartLaplace}
\end{figure}

The line chart (a) in Fig. \ref{songfig:line-chartLaplace} shows the estimation errors over $Pts$ versus the number of samples in different approximation methods for 2-dimensional Laplace equations. While the WoS-driven NN smoothed the discrete approximation results of WoS, our WoS-NN model performed even better than both of them. In WoS-NN, the simulations of local solutions and gradients are integrated and simultaneous, and they optimize each other, enabling superior simulation results. The plot shows that WoS-NN outperformed the conventional WoS with no more than $5\%$ sampled paths and achieved comparable accuracy (average errors around 0.01) to the WoS-driven NN model with only $10\%$ of paths used. When using a similar amount of samples (around 100000 paths), WoS-NN reduces $92\%$ error of WoS, and reduced the error by $52\%$ compared to the WoS-driven NN. Similar tendencies and model features can also be observed in Poisson cases. Line chart (b) in Fig. \ref{songfig:line-chartLaplace} shows the convergence rate of the Relative Root Mean Square Error (RRMSE) of WoS-NN over $Pts$ versus the number of training epochs for the 2-dimensional Laplace equations, where the input training data is $40000$ WoS sampled paths, batch size 1024. The RRMSE converges with the growth of training epochs.

\subsection{2-dimensional Poisson Equations}
For the 2-dimensional Poisson experiments, we tested the following PDE: 
\begin{equation*}
        \Delta u  = 2x \quad \text{on} \quad \Omega=(-1, 1)^2\text{, and } \quad
        u = xy^2, \quad\text{on}\quad \partial\Omega.
\end{equation*}
\vspace{-8mm}
\begin{table}[ht!]
    \centering
    \begin{tabular}{ |>{\centering\arraybackslash}m{2em}|>{\centering\arraybackslash}m{9em}|>{\centering\arraybackslash}m{6em}|>{\centering\arraybackslash}m{5em}|>{\centering\arraybackslash}m{5em}|>{\centering\arraybackslash}m{5em}| } 
 \hline
  &Method & Total Valid Path Amount & Execution Time & Training Error & Average Error on $Pts$ 
   \\ 
 \hline
 1 & Original WoS & 477698 & 74.685 &  & 0.03262\\
 \hline
 2 & WoS-driven NN & $477698$ & 74.685 1.4405 2.4665  & 0.0019 & 0.00893\\
 \hline
 3 & WoS-NN with uniform starting places & 36586 & 30.9249 45.395 
2.5819 & 0.01478 & \makecell{0.008757 \\[-0.8ex] \cmidrule(l{-7pt}r{-7pt}){1-1} 0.05948}\\
 \hline
    \end{tabular}
    \caption{Experimental results for Poisson equation, with the same layout as Table \ref{songtab:Laplace}}
    \label{songtab:Poisson}
\end{table}

 \begin{figure}[ht!]
    \begin{subfigure}{0.44\textwidth}
    \centering
        \includegraphics[width=1\textwidth]{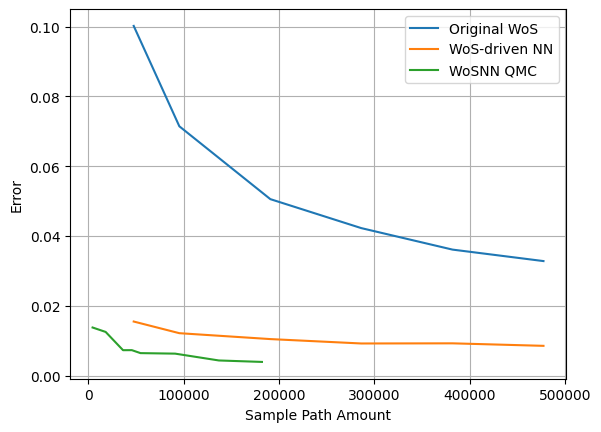}
    \caption{Mean error versus sample path amount}
    \end{subfigure}
    \quad
    \begin{subfigure}{0.42\textwidth}
    \centering
        \includegraphics[width=1\textwidth]{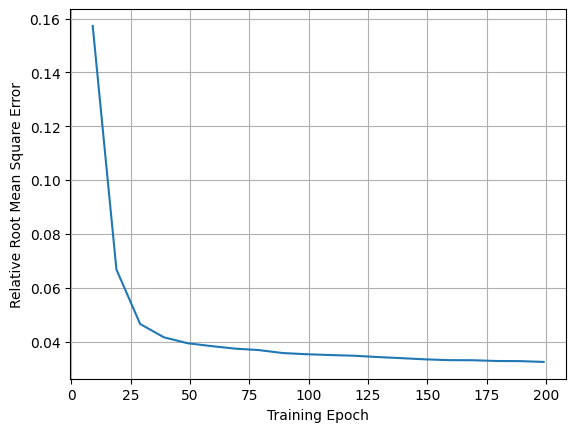}
    \caption{RRMSE versus the training epochs}
    \end{subfigure}
    \caption{Analytical experiment results for the Poisson experiments.}\label{songfig:line-chartPoisson}
\end{figure}
Table \ref{songtab:Poisson} compares representative experiment results on the Poisson equations for different models. All experimental conditions are the same as for the Laplace experiment. Additionally, we used inverse transform sampling for the in-ball sampling in WoS-NN with Green's function density. A probability table is adopted with precision $10^{-5}$. Comparing experiments 1 and 3, our WoS-NN model reduced the error of WoS by $76.23\%$ using only $8\%$ sample paths. Also, as the problem becomes more complex and has higher dimensionality, the time required for sampling begins to dominate over the time needed for training, which is also observed in 3-dimensional experiments. 

\begin{figure}[ht!]
    \begin{subfigure}{0.4\textwidth}
    \centering
        \includegraphics[width=0.85\textwidth]{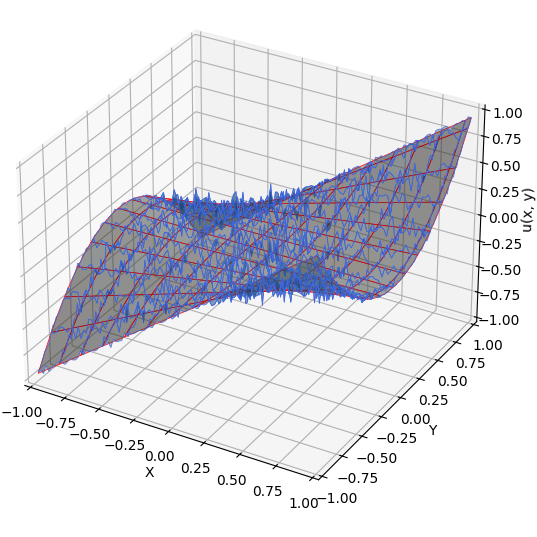}
        \caption{The original WoS method approximation on 9801 points, max = 30 paths per point}
    \end{subfigure}
    \quad
    \begin{subfigure}{0.4\textwidth}
    \centering
        \includegraphics[width=0.85\textwidth]{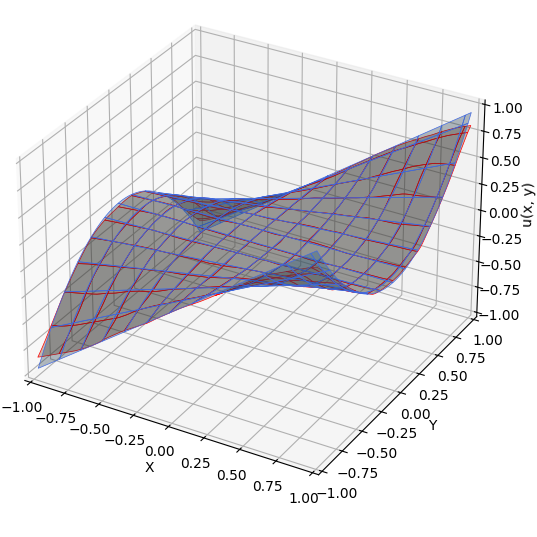}
        \caption{The WoS-driven neural network approximation with 9801 WoS resulting inputs}
    \end{subfigure}

     \begin{subfigure}{0.4\textwidth}
    \centering
        \includegraphics[width=0.85\textwidth]{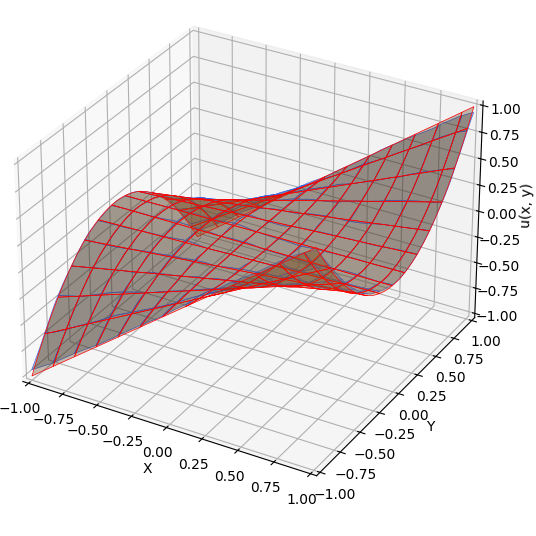}
        \caption{The WoS-NN solution approximation with around 60000 paths sampled by QMC}
    \end{subfigure}
    \quad
     \begin{subfigure}{0.4\textwidth}
    \centering
        \includegraphics[width=0.85\textwidth]{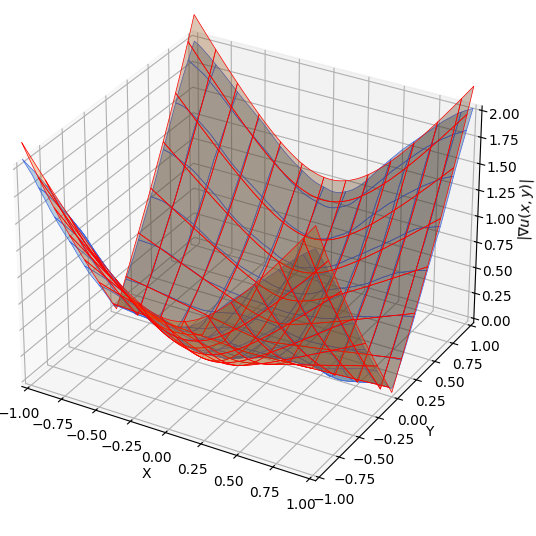}
        \caption{The WoS-NN gradient approximation with around 60000 paths sampled by QMC}
    \end{subfigure}
    \caption{Approximation results with different estimators for 2-D Poisson.\\}
    \label{songfig:2d-Poisson}
    
\end{figure}

Line plot Fig. \ref{songfig:line-chartPoisson} shows the relationship between the sample path amount and approximation error, and the convergence rate of RRMSE on $Pts$ with respect to training epochs for the Poisson experiment on $Pts$ with $40000$ sample paths. The plots indicate that WoS-NN maintained the same superiority in Poisson cases compared to WoS and WoS-driven NN. For experiments with around $100000$ paths, the average Euler distance error of WoS-NN is only $8.82\%$ of the original WoS approximation error and $51.82\%$ of the WoS-driven NN error. To give an intuitive impression, we plotted the Poisson approximation results on $\Omega$ as field approximation, as in Fig. \ref{songfig:2d-Poisson}. For layouts, the red plots present the target function we are approximating, and the blue plots are the approximating results of the current solver. We used gradient norms to exhibit the gradient approximations in the future figures, since gradients are multi-dimensional vectors and are hard to include in a single plot.

\begin{figure}[ht!]
    \begin{subfigure}{0.3\textwidth}
    \centering
        \includegraphics[width=\textwidth]{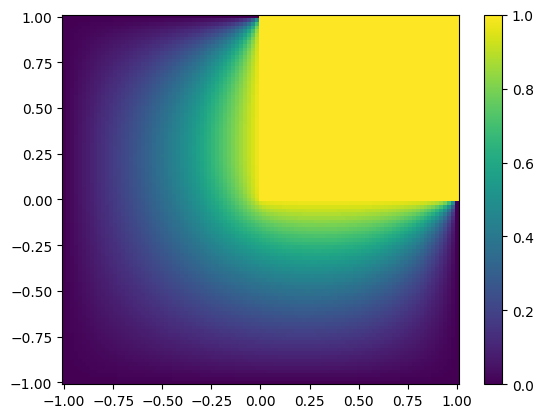}
    \caption{FDM for \eqref{songeqn:L-shape}.}
    \label{songfig:FDML-shaped}
    \end{subfigure}
    \begin{subfigure}{0.3\textwidth}
    \centering
        \includegraphics[width=\textwidth]{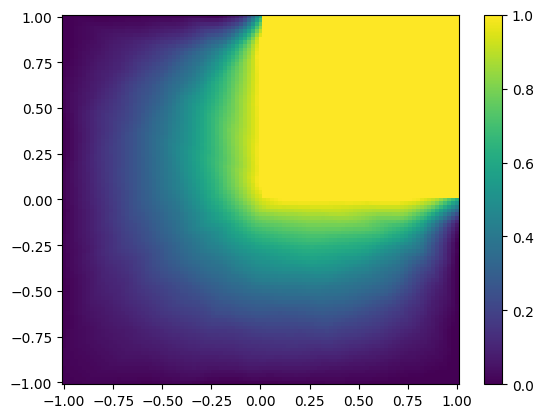}
    \caption{WoS-NN for \eqref{songeqn:L-shape}.}\label{songfig:WoSNNL-shaped}
    \end{subfigure}
    \begin{subfigure}{0.31\textwidth}
    \centering
        \includegraphics[width=\textwidth]{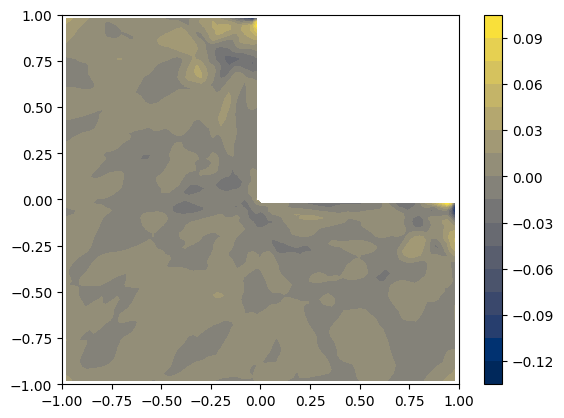} 
    \caption{Error distribution}\label{songfig:DiffL-shaped}
    \end{subfigure}
    \caption{WoS-NN on irregular region and error distribution compared to FDM results} 
\end{figure}

\subsection{Other Experiments}

To test the stability and resilience of our method, we ran WoS-NN on the following PDE, which has broken boundary conditions on an irregular L-shaped region:
\begin{equation}\label{songeqn:L-shape}
    \Delta u = 0 \quad \text{on} \quad \Omega=[-1, 1]^2 \setminus (0, 1]^2,  \text{ and }\quad u = 1\!\!1_{\{x\ge0,y\ge0\}}\quad \text{on} \quad \partial\Omega.
\end{equation}
Here, the Laplace equation has no known closed-form solution due to the irregular shape of the region. The approximation result from the WoS-NN method is as Fig. \ref{songfig:WoSNNL-shaped}. To validate our results, we have experimented with the Finite Difference Method on the same environment as Fig. \ref{songfig:FDML-shaped}. The error distribution of WoS-NN on this L-shaped region, referring to the FDM result, is given as Fig. \ref{songfig:DiffL-shaped}. The mean squared error over the $0.02 * 0.02$ intense grid on the valid region is $MSE = 0.000235$.

\begin{figure}[h]
     \begin{subfigure}{0.4\textwidth}
    \centering
        \includegraphics[width=0.85\textwidth]{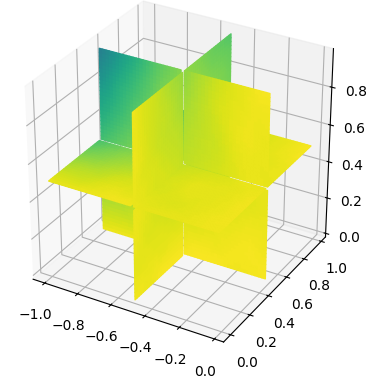}
        \caption{Solution estimation for Poisson equation with source term $2yz$ and boundary condition $x^2yz$}
    \end{subfigure}
    \quad
     \begin{subfigure}{0.4\textwidth}
    \centering
        \includegraphics[width=0.85\textwidth]{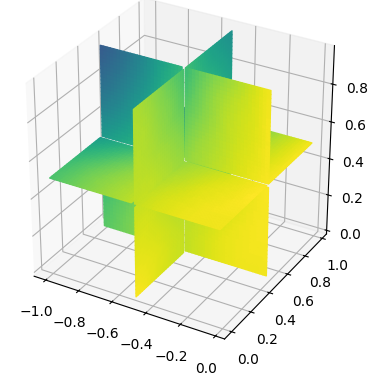}
        \caption{Gradient estimation in norm for Poisson equation with source term $2yz$ and boundary condition $x^2yz$}
    \end{subfigure}
    \caption{WoS-NN result for Laplace and Poisson 3-dimensional on $[-1, 0] \times [0, 1]^2$.}
    \label{songfig:3d-laplace-poisson}
\end{figure}

We also tested WoS-NN on 3-dimensional Poisson equations with $\Omega=[-1,1]^3$, and the result is sketched in Fig. \ref{songfig:3d-laplace-poisson}. In 3-dimensional experiments, WoS paths were sampled with a maximum step of 80 and $\epsilon-$shell of 0.01. 60000 paths were used for training. $YZNN$ is a fully connected feed-forward neural network with three hidden layers $(64, 128, 128)$ and was trained for 50 epochs with batch size 1024, learning rate $2*10^{-4}$, ReLU activation function, and Adam optimizer. To show explicit approximations, we present the results on the region $[-1, 0] \times [0, 1]^2$, where approximations are distinct from each other. We draw middle intersecting planes for both PDEs in the above cube region and plot approximation results on the planes as in Fig. \ref{songfig:3d-laplace-poisson}. 

\section{Conclusion}
In this paper, a novel stochastic solver, WoS-NN, is proposed for solving elliptic PDEs. The method is based on spatial discretization in static PDE regions, using WoS as the sampling tool and neural networks to approximate the PDE solution and gradient. Our test results on 2- and 3-dimensional Laplace and Poisson equations indicate significant improvements in WoS-NN compared to the original WoS using less than $10\%$ of sampling paths. With a similar amount of samples (around 100000 paths), WoS-NN reduces errors by $92\%$ in Laplace and Poisson tests. WoS-NN offers accurate and fast field estimations on both solutions and solution gradients for elliptic PDEs.

\subsection{Future Work}
To further improve and implement our method, several different directions can be considered. First, different parameters and factors in the WoS-NN system can be further optimized, such as WoS step size, training hyperparameters, etc. Second, as a Monte Carlo method, WoS is even more advantageous in high-dimensional scenarios, as well as neural networks. With appropriate samplings and Green's functions, WoS-NN is supposed to be an effective stochastic solver for high-dimensional PDEs. Also, relying on Ito's process and existing extensions of the WoS method \cite{Elepov1969,Rohan2,Simonov2004}, WoS-NN is expected to be scalable to various elliptic problems like Poisson-Boltzmann equations. Finally, we aim to integrate WoS-NN into computer graphics as a robust geometry processing tool \cite{Rohan1}.


\begin{thebibliography}{26}

\bibitem{Beck2019}
Beck, C., E, W., Jentzen, A.: Machine learning approximation algorithms for high-dimensional fully nonlinear partial differential equations and second-order backward stochastic differential equations. Journal of Nonlinear Science 29(4), 1563–1619 (2019). \url{doi:10.1007/s00332-018-9525-3}

\bibitem{Bossy2010}
Bossy, M., Champagnat, N., Maire, S., Talay, D.: Probabilistic interpretation and random walk on spheres algorithms for the Poisson-Boltzmann equation in molecular dynamics. ESAIM: Mathematical Modelling and Numerical Analysis 44(5), 997–1048 (2010). \url{doi:10.1051/m2an/2010050}

\bibitem{Deaconu2013}
Deaconu, M., Herrmann, S.: Hitting time for Bessel processes—walk on moving spheres algorithm (WOMS). The Annals of Applied Probability 23(6), 2259–2289 (2013). \url{http://www.jstor.org/stable/42919710}

\bibitem{dynkin1965}
Dynkin, E., Greenberg, V., Fabius, J., Maitra, A., Majone, G.: Markov Processes: Volume II. Grundlehren der mathematischen Wissenschaften. Springer Berlin Heidelberg (1965). \url{doi:10.1007/978-3-662-00031-1}

\bibitem{Weinan2017}
E, W., Han, J., Jentzen, A.: Deep learning-based numerical methods for high-dimensional parabolic partial differential equations and backward stochastic differential equations. Communications in Mathematics and Statistics 5(4), 349–380 (2017). \url{doi:10.1007/s40304-017-0117-6}

\bibitem{Elepov1969}
Elepov, B., Mikhailov, G.: Solution of the Dirichlet problem for the equation $\Delta u-cu = -q$ by a model of “walk on spheres.” USSR Computational Mathematics and Mathematical Physics 9(3), 194–204 (1969). \url{doi:10.1016/0041-5553(69)90070-6}

\bibitem{Ito}
Itô, K.: On a formula concerning stochastic differentials. Nagoya Mathematical Journal 3, 55–65 (1951). \url{doi:10.1017/S0027763000012216}

\bibitem{Kakutani1944}
Kakutani, S.: Two-dimensional Brownian motion and harmonic functions. Proceedings of the Imperial Academy 20(10), 706 – 714 (1944).  \url{doi:10.3792/pia/1195572706}

\bibitem{Siggraph1}
Li, Z., Yang, G., Deng, X., De Sa, C., Hariharan, B., Marschner, S.: Neural caches for Monte Carlo partial differential equation solvers. In: SIGGRAPH Asia 2023 Conference Papers, SA’23. Association for Computing Machinery, New York, NY, USA (2023). \url{doi:10.1145/3610548.3618141}

\bibitem{Siggraph2}
Li, Z., Yang, G., Zhao, Q., Deng, X., Guibas, L., Hariharan, B., Wetzstein, G.: Neural control variates with automatic integration. In: ACM SIGGRAPH 2024 Conference Papers, SIGGRAPH ’24. Association for Computing Machinery, New York, NY, USA (2024). \url{doi:/10.1145/3641519.3657395}

\bibitem{Simonov2004}
Mascagni, M., Simonov, N.A.: Monte Carlo methods for calculating some physical properties of large molecules. SIAM Journal on Scientific Computing 26(1), 339–357 (2004). \url{doi:10.1137/S1064827503422221}

\bibitem{Muller1956}
Muller, M.E.: Some Continuous Monte Carlo Methods for the Dirichlet Problem. The Annals of Mathematical Statistics 27(3), 569 – 589 (1956). \url{doi:10.1214/aoms/1177728169}

\bibitem{Capacitance1991}
Nabors, K., White, J.: Fastcap: a multipole accelerated 3-d capacitance extraction program. IEEE Transactions on Computer-Aided Design of Integrated Circuits and Systems 10(11), 1447–1459 (1991). \url{doi:10.1109/43.97624}

\bibitem{Nam2024}
Nam, H.C., Berner, J., Anandkumar, A.: Solving Poisson equations using neural walk-on-Spheres. In: Proceedings of the 41st International Conference on Machine Learning, ICML’24. JMLR.org (2024)

\bibitem{numeric}
Pinder, G.F.: Numerical Methods for Solving Partial Differential Equations: A Comprehensive Introduction for Scientists and Engineers. John Wiley and Sons, Inc.: Wiley (2018)

\bibitem{RAISSI2019}
Raissi, M., Perdikaris, P., Karniadakis, G.: Physics-informed neural networks: A deep learning framework for solving forward and inverse problems involving nonlinear partial differential equations. Journal of Computational Physics 378, 686–707 (2019). \url{doi:10.1016/j.jcp.2018.10.045}

\bibitem{Sabelfeld1991}
Sabelfeld, K.K.: Monte Carlo methods in boundary value problems (1991). \url{https://api.semanticscholar.org/CorpusID:117752869}

\bibitem{Sabelfeld2013}
Sabelfeld, K.K., Shalimova, I.A.: Spherical and Plane Integral Operators for PDEs. De Gruyter, Berlin, Boston (2013). \url{doi:10.1515/9783110315332}

\bibitem{SalsaPDE}
Salsa, S.: Partial Differential Equations in Action: From Modelling to Theory. Springer Milano(2008). \url{doi:10.1007/978-88-470-0752-9}

\bibitem{Rohan1}
Sawhney, R., Crane, K.: Monte Carlo geometry processing: a grid-free approach to PDE-based methods on volumetric domains. ACM Transactions on Graphics 39(4) (2020). \url{doi:10.1145/3386569.3392374}

\bibitem{Rohan2}
Sawhney, R., Seyb, D., Jarosz, W., Crane, K.: Grid-free Monte Carlo for PDEs with spatially varying coefficients (2022). \url{doi:10.1145/3528223.3530134}

\bibitem{Simonov2007}
Simonov, N.A., Mascagni, M., Fenley, M.O.: Monte Carlo-based linear Poisson-Boltzmann approach makes accurate salt-dependent solvation free energy predictions possible. The Journal of Chemical Physics, 127(18), 185105 (2007). \url{doi.org/10.1063/1.2803189}

\bibitem{Sirignano2017}
Sirignano, J., Spiliopoulos, K.: Dgm: A deep learning algorithm for solving partial differential equations. Journal of Computational Physics 375, 1339–1364 (2018). \url{doi:10.1016/j.jcp.2018.08.029}


\end{thebibliography}
%

\end{document}